\documentclass[10pt]{amsart}
\usepackage{verbatim}
\usepackage{eucal}
\usepackage{amssymb}
\usepackage{mathrsfs}
\usepackage{graphicx}
\usepackage{psfrag}
\usepackage{xypic}
\newif \ifwide
\newif \ifavnermargin
\def \makemargins{
\ifwide
        \oddsidemargin .25in
        \evensidemargin .25in
        \textwidth 6.00in
\else
\fi
\ifavnermargin
        \headheight=7pt
        \textheight=574pt
        \textwidth=432pt
        \topmargin=14pt
        \oddsidemargin=18pt
        \evensidemargin=18pt
\else   
\fi
}

\avnermargintrue
\makemargins

\theoremstyle{plain}
\newtheorem{theorem}[subsection]{Theorem}
\newtheorem{proposition}[subsection]{Proposition}
\newtheorem{lemma}[subsection]{Lemma}

\theoremstyle{definition}
\newtheorem{definition}[subsection]{Definition}

\theoremstyle{remark}
\newtheorem{remark}[subsection]{Remark}

\newcount\timehh\newcount\timemm
\timehh=\time 
\divide\timehh by 60 \timemm=\time
\newcommand{\draftauthor}[1]{\author{#1
    {
      --- \protect \protect\sc\today\ ---
      \ifnum\timehh<10 0\fi\number\timehh\,:\,\ifnum\timemm<10 0\fi\number\timemm
      \protect \, \, \protect \bf DRAFT
    }
  }
}
\newcommand{\PP}{{\mathbb P}}
\newcommand{\CC}{{\mathbb C}}
\newcommand{\ZZ}{{\mathbb Z}}
\newcommand{\QQ}{{\mathbb Q}}

\newcommand{\ee}{{\mathrm{e}}}
\newcommand{\ii}{{\mathrm{i}}}

\newcommand{\HHH}{{\mathfrak{H}}}
\newcommand{\TTT}{{\mathscr{T}}}

\newcommand{\EEE}{{\mathscr{E}}}
\newcommand{\MMM}{{\mathscr{M}}}
\newcommand{\SSS}{{\mathscr{S}}}
\newcommand{\ac}[1]{{\rm a.c.}\Bigl(#1\Bigr)}
\newcommand{\stroke}[2]{{#1\!\bigm|\!{#2}}}

\renewcommand{\theta}{\vartheta}
\renewcommand{\tilde}{\widetilde}

\renewcommand{\mod}{\bmod}

\DeclareMathOperator{\Hom}{Hom}

\DeclareMathOperator{\codim}{codim}

\begin{document}

\title{Toric modular forms and nonvanishing of $L$-functions}

\newif \ifdraft
\def \makeauthor{
\ifdraft
        \draftauthor{Lev A. Borisov and Paul E. Gunnells}
\else
\author{Lev A. Borisov}
\address{Department of Mathematics\\
Columbia University\\
New York, NY  10027}
\email{lborisov@math.columbia.edu}

\author{Paul E. Gunnells}
\address{Department of Mathematics and Computer Science\\
Rutgers University\\
Newark, NJ  07102}
\email{gunnells@andromeda.rutgers.edu}
\fi
}

\thanks{Both authors were partially supported by a Columbia
University Faculty Research grant}

\draftfalse
\makeauthor

\ifdraft
        \date{\today}
\else
        \date{October 26, 1999. Revised September 15, 2000}
\fi

\subjclass{11F11, 11F25, 14M25}
\keywords{Toric varieties, modular forms, Manin symbols, special
values of $L$-functions}

%
%

\begin{abstract}
In a previous paper \cite{BorGunn}, we defined the space of toric
forms $\TTT(l)$, and showed that it is a finitely generated subring
of the holomorphic modular forms of integral weight on the congruence
group $\Gamma_1(l)$.  In this article we prove the following theorem:
modulo Eisenstein series, the weight two toric forms coincide exactly
with the vector space generated by all cusp eigenforms $f$ such that
$L(f,1) \not = 0$.  The proof uses work of Merel, and involves an
explicit computation of the intersection pairing on Manin symbols.
\end{abstract}
\maketitle

%
%
\section{Introduction}\label{introduction}
\subsection{}
Let $l>1$ be an integer, let $\MMM (l)$ be the space of weight two
modular forms on the congruence group $\Gamma _{1} (l) \subset SL_{2}
(\ZZ )$, and let $\SSS (l)$ be the subspace of cusp forms.  Let $f\in
\SSS (l)$ be an eigenform for the Hecke operators $T_{p}$, where $p$
is coprime to $l$, and let $L (f,s)$ be the associated $L$-function.
Then the order of vanishing of $L (f,s)$ at $s=1$ is called the
\emph{analytic rank} of $f$.  This terminology comes from the Birch
and Swinnerton-Dyer conjecture, which asserts that the analytic rank
times $\dim A_{f}$ is the same as rank of the group $A_{f} (\QQ)$,
where $A_{f}$ is the abelian variety associated to $f$ by the
Eichler-Shimura construction \cite{Knapp, Tate}.

\subsection{}
In this paper we present an elementary construction of the subspace of
$\SSS (l)$ spanned by forms of analytic rank zero.  Our main result
(Theorem~\ref{main}) is that, modulo Eisenstein series, this space is
isomorphic as a Hecke module to the space $\TTT_2 (l)$ of weight two
\emph{toric modular forms} of level $l$.  These modular forms were
constructed and studied in \cite{BorGunn}, where we presented explicit
generators of $\TTT (l)$ and described their $q$-expansions at
infinity.  The construction of $\TTT (l)$ and its relevant properties
are summarized in Theorem~\ref{summary}.  For the remainder of
this introduction, we describe the proof of Theorem~\ref{main}.

\subsection{}
First, in \S\ref{manin.symbols} we recall results about the
\emph{Manin symbols}.  We discuss various homology groups associated
to the modular curve in terms of modular symbols, and describe the
intersection pairing.  We define the space of \emph{plus}
(respectively \emph{minus}) symbols $M_{+}$ (resp. $M_{-}$), and their
\emph{cuspidal} subspaces $S_{+}$ and $S_{-}$, and we describe the
intersection pairing in terms of Manin symbols (Proposition \ref{2.9}).
We finish this
section by recalling Merel's description of the Hecke action on Manin
symbols (Theorem~\ref{HeckeisHecke}).

Next, in \S\ref{toric.forms.section} we recall results about toric
modular forms from \cite{BorGunn}, specialized to the case of weight
two.  Using the toric forms, we define a map $\mu \colon
M_{-}\rightarrow \MMM (l)/\EEE (l)$, where $\EEE (l)\subset \MMM (l)$
is the subspace of Eisenstein series.  The main result
(Theorem~\ref{muisHecke}) is that $\mu $ is Hecke-equivariant.  We
also describe the composition of $\mu $ with the \emph{Fricke
involution} $W_{l}$.

Finally, in \S\ref{rank.zero.section} we define a map $\rho \colon \SSS
(l)\rightarrow \SSS (l)$ whose image is spanned by Hecke eigenforms of
analytic rank zero.  Then we put all these maps together
to get a sequence
\[
\xymatrix@1{{\SSS (l)}\ar[r]^{\int }&M_{+}^{*}\ar[r]^{\pi
}&M_{-}\ar[r]^-{W_{l}\circ \mu }&\MMM (l)/\EEE (l)\ar[r]^-{\sim}&\SSS (l)},
\]
where the final map is the Hecke-equivariant isomorphism between $\MMM
(l)/\EEE (l)$ and $\SSS (l)$.  Then in Theorem~\ref{keytheorem} we
show that this sequence equals $\rho $, from which we obtain
Theorem~\ref{main}.

The Eisenstein series $s_{a/l}$ have appeared in the literature before.
They were originally studied by Hecke \cite{hecke,Lang} and have
recently appeared in the work of Kato on Euler systems
\cite{kato.lectures,Scholl}. Moreover, it might be possible
to derive Theorem \ref{main} using the formulas of 
\cite[\S 4]{Scholl}, which were obtained 
using the Rankin-Selberg method.

\subsection{Acknowledgments}
We thank Dorian Goldfeld for encouraging this work, and Siman Wong for
some interesting discussions. We thank Kazuya Kato and Anthony Scholl
for helpful correspondence. Last, but not least, we thank the referee
for suggesting a possible connection of our work to
\cite{kato.lectures,Scholl}, for providing revisions that
substantially improved the exposition, and for the statement and proof
of Proposition \ref{ref.prop}.  After this paper was submitted, Lo\"\i
c Merel informed us that \cite{merel2} contains an intersection
formula for $\Gamma_0(l)$ analogous to our formula of Proposition
\ref{2.9}.

\section{Manin symbols}
\label{manin.symbols}

\subsection{}
Let $l>1$ be an integer, and let $\Gamma _{1} (l)\subset SL_{2} (\ZZ
)$ be the subgroup of matrices congruent to
$\left(\begin{smallmatrix}1&*\\0&1\end{smallmatrix}\right)$ mod $l$.
Let $\HHH $ be the upper halfplane, and let $\HHH ^{*} = \HHH \cup \PP
^{1} (\QQ)$ be the usual partial compactification obtained by
adjoining cusps.  Let $X _{1}(l) = \Gamma _{1} (l)\backslash\HHH ^{*}$
be the modular curve with cusps $\partial X _{1}(l) $, and let $Y_{1}
(l) := X _{1}(l)\smallsetminus \partial X _{1}(l) $.

Let $M$ be the relative homology $H_{1} (X _{1}(l), \partial X
_{1}(l); \CC)$, and let $S\subset M$ be the subspace $H_{1} (X_{1}
(l), \CC )$.  The intersection product induces a perfect pairing of
complex vector spaces
$$
H_1(X_1(l),\partial X_1(l);\CC)\times H_1(Y_1(l),\CC)\longrightarrow
\CC,
$$
which allows us to identify $M^{*} = \Hom _{\CC } (M,\CC )$ with
$H_1(Y_1(l),\CC)$.

\subsection{}
Manin's theory \cite{Manin} gives a concrete description of $M$ as
follows.  For any $\alpha, \beta$ in $\PP^1(\QQ)$, let $\{\alpha
,\beta  \}\in M$ be the class of the image of a continuous path in
$\HHH ^{*}$ from $\alpha $ to $\beta $.  Then the $\CC$-linear map 
\begin{align*}
\CC[\pm\Gamma_1(l)\backslash SL_2(\ZZ)]&\longrightarrow M\\
\pm\Gamma_1(l)g &\longmapsto \{g0,g\infty\}
\end{align*}
is well defined and surjective. If we denote the basis element
corresponding to the coset $x$ by 
$[x]$, then the kernel of this map is generated by
elements of the form  $[x]+[x\sigma]$ and
$[x]+[x\tau]+[x\tau^2]$. Here $x$ runs through
$\pm\Gamma_1(l)\backslash SL_2(\ZZ)$, and $\sigma, \tau$ are
elements of $SL_2(\ZZ)$ (of order $4$ and $3$) that stabilize $\ii$ and
$\rho=\ee^{2\pi\ii/3}$ respectively, and satisfy $\sigma0=\infty$
and $\tau\infty=0$.

By duality one can describe $M^*$ as the subspace of
$\CC[\pm\Gamma_1(l)\backslash SL_2(\ZZ)]$ generated by elements
$\sum_x\lambda_x[x]$ satisfying $\lambda_x+\lambda_{x\sigma}=0$ and
$\lambda_x+\lambda_{x\tau}+\lambda_{x\tau^2}=0$.  One can realize this
description geometrically as follows. Consider the geodesic path in
the upper half-plane from $\ii$ to $\rho$. For $g\in SL_2(\ZZ)$, the
image in $Y_1(l)$ of its translates by $g$ depends only on the coset
$\pm \Gamma_1(l)g$.  Denote by $c_{\Gamma_1(l)g}$ the arc in $Y_1(l)$
associated to this image. One can show easily that
$\sum_x\lambda_xc_x$ is a closed cycle if and only if
$\lambda_x+\lambda_{x\sigma}=0$ and
$\lambda_x+\lambda_{x\tau}+\lambda_{x\tau^2}=0$. In that case,
$\sum_x\lambda_xc_x$ belongs to $M^*=H_1(Y_1(l),\CC)$. Using easy
considerations on fundamental domains, one shows that the intersection
pairing of the image in $M$ of $[y]$ and of $\sum_x\lambda_xc_x$ is
equal to $\lambda_y$. For details, we refer to \cite{merel3}.

\subsection{}\label{pi.explicit}
We will now calculate explicitly the natural map 
$\pi\colon M^*\to M$, 
which is the composition of canonical topological maps
$$
M^*=H_1(Y_1(l),\CC)\rightarrow H_1(X_1(l),\CC)
\rightarrow H_1(X_1(l),\partial X_1(l);\CC)=M.
$$

\begin{proposition}\label{ref.prop}
Let $x\mapsto g_x$ be a section of the surjective map
$SL_2(\ZZ)\rightarrow\pm \Gamma_1(l)\backslash SL_2(\ZZ)$.  
Let $\sum_x\lambda_xc_x\in M^*$.
The map $\pi$ is given
by the following formulas:
$$
\pi(\sum_x\lambda_xc_x)={1\over
6}\sum_x(\lambda_{x\tau}-\lambda_{x\tau^2})\{g_x0,g_x\infty\}
=\frac 16 \sum_x(\lambda_{x\sigma\tau\sigma}-\lambda_{x\tau^2})
\{g_x0,g_x\infty\},$$
where the sums are taken over $x\in \pm \Gamma_1(l)\backslash SL_2(\ZZ)$.
\end{proposition}

\begin{proof}
This is a simple computation.
One has
$$\pi(\sum_x\lambda_xc_x)=\sum_x \lambda_x\{g_x\ii,g_x\rho\}
=\sum_x\lambda_x\{g_x\ii,g_x\infty\}-\sum_x\lambda_x\{g_x\rho,g_x\infty\}.
$$
Here we are abusing the notation $\{\alpha,\beta\}$, which now
also denotes the arc in $X_1(l)$ that is the image of an arc in $\HHH^*$.

We use the fact that $g_{x\sigma}$ and $g_{x}\sigma$ lie in the same
coset of $\pm\Gamma_1(l)$ (and similarly for $\tau$) to rewrite the
right hand side of the above expression as 
$$
{1\over 2}\sum_x(\lambda_x\{g_x\ii,g_x\infty\}
+\lambda_{x\sigma}\{g_x\sigma \ii,g_x\sigma\infty\})
-{1\over 3}\sum_x(\lambda_x\{g_x\rho,g_x\infty\}
+\lambda_{x\tau}\{g_x\tau\rho,g_x\tau\infty\}
+\lambda_{x\tau^2}\{g_x\tau^2\rho,g_x\tau^2\infty\}).
$$
Using the relations $\sigma \ii=\ii$,
$\tau \rho=\rho$, $\lambda_x+\lambda_{x\sigma}=0$ and
$\lambda_x+\lambda_{x\tau}+\lambda_{x\tau^2}=0$, we have
$$
{1\over 2}\sum_x(\lambda_x\{g_x\infty,g_x\infty\}
+\lambda_{x\sigma}\{g_x\infty,g_x\sigma\infty\})
-{1\over 3}\sum_x(\lambda_x\{g_x\infty,g_x\infty\}
+\lambda_{x\tau}\{g_x\infty,g_x\tau\infty\}
+\lambda_{x\tau^2}\{g_x\infty,g_x\tau^2\infty\}).
$$
Using the relations $\tau\infty=0$ and $\{g_x\infty,g_x\infty\}=0$ and
reindexing the last set of terms of the sum, we obtain
$$
{1\over2}\sum_x\lambda_x\{g_x0,g_x\infty\}
-{1\over 3}\sum_x\lambda_{x\tau}\{g_x\infty,g_x0\}
-{1\over 3}\sum_x\lambda_{x}\{g_x\tau\infty,g_x\infty\}.
$$
Using again the relations $\tau\infty=0$,
$\{g_x\infty,g_x0\}=-\{g_x0,g_x\infty\}$, and
$\lambda_x+\lambda_{x\tau}+\lambda_{x\tau^2}=0$, we arrive at the
first equality of the statement. To get the second equality, one
changes the index to $x\sigma$ and uses the relation
$\lambda_x+\lambda_{x\sigma}=0$.
\end{proof}

\begin{remark}
Proposition \ref{ref.prop} remains valid if 
$\Gamma_1(l)$ is replaced by any finite-index subgroup of $SL_2(\ZZ)$.
\end{remark}

\subsection{}
The elements of  $\Gamma_1(l)\backslash SL_2(\ZZ)$ can be
represented 
by pairs $(u,v)$, where $u,v\in \ZZ/l\ZZ$ and $g.c.d.(u,v,l)=1$.
The discussion above
shows that $M$ can be described as the 
$\CC$-vector space generated by the symbols $(u,v)$ modulo the relations
\begin{enumerate}
\item $(u,v) + (-v,u) = 0$.
\item $(u,v) + (v,-u-v) + (-u-v,v) = 0$.
\end{enumerate}
Pairs $(u,v)$ are called \emph{Manin symbols}.  
Two subspaces of $M$ will play an important role in what follows.  Let
$\iota \colon M \rightarrow M$ be the involution that takes
$(u,v)\mapsto (-u,v)$.  
\begin{definition}\label{plusnminus}
The space of \emph{plus symbols} $M_{+}\subset M$ is the subspace
consisting of symbols $x$ satisfying $\iota (x) = x$.  Similarly,
the space of \emph{minus symbols} $M_{-}\subset M$ is the subspace
consisting of symbols $x$ satisfying $\iota (x) = -x$.
\end{definition}

We have \emph{symmetrization maps} $(\phantom{a},\phantom{a})_{\pm }\colon M
\rightarrow M_{\pm }$ given by $(u,v)_{\pm } := ((u,v)\pm (-u,v))/2$.
We also introduce the corresponding 
spaces of cuspidal symbols $S_\pm\subseteq M_\pm$.
The spaces $S_{\pm }$ can each be seen as the dual of the space of
cusp forms as follows.  Let $\MMM (l)$ be the $\CC $-vector space of
weight two holomorphic modular forms on $\Gamma_{1} (l) $, and let
$\SSS (l)\subset \MMM (l)$ be the subspace of cusp forms.  Let
$(u,v)\in M$, and let the cusps corresponding to $u, v$ be $\alpha ,
\beta $ respectively.  The pair $\alpha , \beta $ induces a geodesic
on $X_{1} (l)$; hence given any $f\in \SSS (l)$, we can form the
integral
\[
\int _{\alpha }^{\beta }f (z)\,dz \in \CC,
\]   
which converges since $f$ is a cusp form.  In this way we identify
Manin symbols with functionals on cusp forms, and likewise cusp forms
with elements of the dual space $M^{*}$.  We obtain a pairing
\begin{align*}
M\times \SSS (l)&\longrightarrow \CC,\\
((u,v), f)&\longmapsto \langle f, (u,v)\rangle.
\end{align*}
In general this pairing is
degenerate, although we have the following result:

\begin{proposition}
\cite[Proposition 8]{Merel}
The pairings 
\[
S_{\pm } \times \SSS (l)\rightarrow \CC 
\]
are nondegenerate.
\end{proposition}

\begin{remark}\label{iotageo}
The involution $(u,v)\mapsto (-u,v)$ on Manin symbols is induced from
the action of the map $\tau \mapsto -\bar \tau , \tau \in \HHH $ on
geodesics.
\end{remark}

\subsection{}
Let $(u,v)$ be a Manin symbol.  For any $\varphi \in
M^{*}$, we define $\varphi $ on ``degenerate'' symbols $(u,v)$ with
$\ZZ u+\ZZ v \not = \ZZ /l\ZZ $ by setting $\varphi (u,v) = 0$.  
This convention is somewhat artificial but turns out to be quite 
useful.

We now rewrite the map $\pi$ of Proposition \ref{ref.prop}
on $(M_+)^*$ in a form that will be useful later.

\begin{proposition}\label{2.9}
The image of $M_+^*$ under $\pi$ is $S_-$. For
any element of $\varphi\in (M_+)^*$, we have
$$
\pi(\varphi)=\frac 1{12}\sum_{a,b=0}^{l-1}
\varphi ((a,a-b)_{+}-(a,a+b)_{+}) (a,b)_{-}.
$$
In addition, $\pi(M_+^*)=S_-$.
\end{proposition}

\begin{proof}
Because $\iota$ comes from the orientation-reversing automorphism 
of $X_1(l)$, it anticommutes with $\pi$. The surjectivity of 
the map $M^*=H_1(Y_1(l),\CC)\rightarrow H_1(X_1(l),\CC)$
then implies $\pi(M_+^*)=S_-$. The second part of the statement follows
from the second equality in Proposition \ref{ref.prop} and the
definitions of the symmetrization maps. The coefficient is changed
to $\frac 1{12}$ because the sum is now over $\Gamma_1(l)\backslash
SL_2(\ZZ)$ instead of $\pm\Gamma_1(l)\backslash
SL_2(\ZZ)$.
\end{proof}

\subsection{}
To conclude this section, we present Merel's description of the Hecke
action on the Manin symbols.  Let $n\geq 1$ be an integer, and let
$T_{n}$ be the associated Hecke operator (cf. \cite{Lang}).  We denote
the action of $T_{n}$ on a modular form $f$ by $\stroke{f}{T_{n}}$.

\begin{theorem}\label{HeckeisHecke}
\cite[Theorem 2 and Proposition 10]{Merel} The operator $T_{n}$ acts on
any Manin symbol
$(u,v)$ via
\begin{equation}\label{heckeact}
T_{n} (u,v) = \sum _{\substack{a>b\geq 0\\
d>c\geq 0\\
ad-bc = n}} (au+cv, bu+dv).
\end{equation}
If $n$ is not coprime to $l$, then we omit the terms for which
$g.c.d.(l,au+cv,bu+dv ) > 1$.
This action is compatible with the pairing between cusp forms and
Manin symbols
\[
\langle \stroke{f}{T_{n}}, (u,v) \rangle= \langle f, T_{n} (u,v)\rangle.
\]
\end{theorem}

It is also easy to show that this Hecke action is compatible with the
symmetrization maps:
\begin{proposition}\label{HeckeonSymbols}
\[
T_{n} ((u,v)_{\pm }) = (T_{n} (u,v))_{\pm }.
\]
\end{proposition}

\begin{proof}
This follows from switching $a$ with $d$ and $b$ with $c$ in
\eqref{heckeact}.
\end{proof}

\section{Toric forms of weight two}\label{toric.forms.section}

\subsection{} 
Let us briefly review the contents of \cite{BorGunn}.  For every
integer $l>1$, we defined a certain Hecke-stable subring of the ring
of modular forms for $\Gamma _{1} (l)$, called the subring of {\it
toric forms} $\TTT (l)$. In the present paper we are only concerned
with weight two toric forms, which greatly simplifies the
combinatorial data needed to encode toric varieties.

Let $N=\ZZ^2$ be a lattice of rank two, where lattice simply means 
a free abelian group. A (compact) toric variety of dimension two is
defined uniquely by a collection of $k$ rational rays from the origin,
such that the angle between any two consecutive rays is less than $\pi$. 
We denote the minimum nonzero lattice points on these rays by $d_i$,
$i=0,\dots ,k$, where $i$ increases counterclockwise, and where $d_{0}
= d_{k}$.

To every such collection one associates a {\it fan} $\Sigma$, which is   
a collection of $2k+1$ cones in $N_\QQ$.   This fan contains $k$
two-dimensional cones 
\[
\{\QQ_{\geq 0}d_i+\QQ_{\geq 0}d_{i+1}\mid i=0,\dots ,k-1\},
\]
$k$ one-dimensional cones 
\[
\{\QQ_{\geq 0}d_i\mid i=0,\dots ,k-1\},
\]
and one zero-dimensional 
cone $\{0\}$.  The corresponding toric variety is 
smooth if and only if $(d_i,d_{i+1})$ is a basis of $\ZZ^2$ for 
every $i$.

To define a toric form we need an additional piece of data, namely a
{\it degree function} with respect to $\Sigma $.  This is a
piecewise-linear function $\deg\colon N\to \QQ$ that is linear on the
cones of $\Sigma$.  Every such function is uniquely determined by the 
values $\alpha_i = \deg ( d_i)$. 

\begin{definition}
\cite{BorGunn} Suppose $\deg \colon N \rightarrow \QQ $ is a degree
function with respect to the fan $\Sigma$, and that $\alpha _{i}\not
\in\ZZ $ for all $i$.  Then the \emph{toric form} associated to
$(N,\deg)$ is the function $f_{N,\deg}\colon \HHH\rightarrow \CC $
defined by
$$
f_{N,\deg}(q) := \sum_{m\in M}\Bigl(\sum_{C\in\Sigma}(-1)^{\codim C}
\ac{\sum_{n\in C} q^{m\cdot n} \ee^{2\pi\ii \deg(n)}}\Bigr).
$$
Here $M={\rm Hom}(M,\ZZ)$ is the dual of $N$, $q=\ee ^{2\pi \ii \tau
}$ where $\tau \in \HHH $, and $\text{a.c.}$ means analytic
continuation of a sum from its region of convergence to all $q$ and $m$.
\label{modularfromdegree}
\end{definition}

It turns out that $f_{N,\deg}$ does not change if $\Sigma$ is subdivided,
and thus does not depend on $\Sigma $.  This is why $\Sigma$ is
omitted from the notation.  A toric form is, by definition, any linear
combination of $f_{N,\deg}$.  We will now state in one theorem most of 
the results of \cite{BorGunn}, specialized to the case of weight two.

\begin{theorem}
\cite{BorGunn} Suppose that $\deg(N)\subseteq \frac1l\ZZ$, and that $\alpha
_{i}\not \in \ZZ $ for all $i$.  Then $f_{N,\deg}$ is a holomorphic
modular form of weight two with respect to $\Gamma_1(l)$. If $l\geq 5$
then it is a linear combination of pairwise products of the forms
$s_{a/l}$, $a = 1, \dots , l-1$, where
\[
s_{a/l}(\tau)=\frac 1{2\pi\ii} \partial_z(\log\theta)(\frac
al,\tau).
\]
Here $\theta(z,\tau)$ is the standard theta
function \cite[Chapter 5]{Chandra}.  If $l<5$, then
the span $\TTT _{2} (l)$ of all toric forms of level $l$ and weight
two coincides with the space of all modular forms of weight two; in
particular, it consists only of Eisenstein series.  The space $\TTT
_{2} (l)$ is stable under the action of Hecke operators and the Fricke
involution.  The space of all weight two toric forms of all levels is
stable under Atkin-Lehner liftings $f(\tau)\mapsto f(n\tau)$.
\label{summary}
\end{theorem}

Let $p$ be a prime not dividing $l$.  We will need an explicit
formula for the action of the Hecke operator $T_p$ on $f_{N, \deg}$.
This follows immediately from a formula in \cite[Theorem 5.3]{BorGunn}
specialized to the case of weight two.
\begin{proposition}
Let $f_{N,\deg}$ be a toric form of weight two. Then 
$$
\stroke{f_{N,\deg}}{T_p} = \sum_{S}f_{S,p\deg},
$$
where the sum is taken over all superlattices $S\supseteq N$ with $[S:N]=p$.
\label{explicitHecke}
\end{proposition}

\subsection{}
It is not hard to write down the explicit $q$-expansions of $s_{a/l}$.
\begin{proposition}
Denote $w=\exp(2\pi\ii/l)$. Then
$$s_{a/l}(q)=\frac{w^a+1}{2(w^a-1)} - 
\sum _{d}q^{d}\sum _{k|d} (w^{ka}-w^{-ka}).
$$
\label{explicitsa}
\end{proposition}

\begin{proof}
It is easy to compute the logarithmic derivative of $\theta$ using the
Jacobi triple product formula \cite[Chapter
5, Theorem 6]{Chandra}. Details are
left to the reader.
\end{proof}

In \cite{BorGunn} we introduced weight two modular forms 
$s_{a/l}^{(2)}$ given by 
$$
s_{a/l}^{(2)}(q)=
-\frac 1{4\pi^2} 
\left(\frac{\partial^2}{\partial z^2}\right)_{z=0}{\rm log}
\left(
{
\frac
{z\theta(z+a/l)\theta'(0)}
{\theta(z)\theta(a/l)}
}
\right)
.
$$
These are also toric forms, and the following relations allow one to
express them as linear combinations of products $s_{a/l} s_{b/l}$ for
$l\geq 5$.
\begin{proposition}
If $a+b+c = 0\mod l$ and $a,b,c \neq 0\mod l$ then 
$$
s_{a/l}s_{b/l}+
s_{b/l}s_{c/l}+
s_{c/l}s_{a/l}
=-\frac 12
\left(
(s_{a/l})^2+
(s_{b/l})^2+
(s_{c/l})^2+
s_{a/l}^{(2)}+
s_{b/l}^{(2)}+
s_{c/l}^{(2)}
\right).
$$
\label{relation}
\end{proposition}

\begin{proof}
This is a consequence of a more general formula of
\cite{BorGunn}.
\end{proof}

\begin{proposition}
The modular forms $(s_{a/l})^2$ and 
$s_{a/l}^{(2)}$ are Eisenstein series for all $a$.
\label{toricEisenstein}
\end{proposition}

\begin{proof}
As in Proposition \ref{explicitsa}, we get 
$$(2\pi\ii)^{-1}\partial_z \log\theta(z,\tau)=
\frac{\ee^{2\pi\ii z} +1}{2(\ee^{2\pi\ii z}-1)}
-
\sum_{d>0}q^d\sum_{k|d}
(\ee^{2\pi\ii kz}-\ee^{-2\pi\ii kz}).$$
Differentiating it again with respect to
$z$ and plugging in $z=a/l$, we get
$$
(2\pi\ii)^{-2}\left(\frac {\partial^2}{\partial z^2}\right)
\log\theta(\frac al,\tau)
=
\frac{\ee^{2\pi\ii a/l}}
{(\ee^{2\pi\ii a/l}-1)^2}
-\sum_{d>0}q^d\sum_{k|d}k(\ee^{2\pi\ii ka/l}+\ee^{-2\pi\ii ka/l}).
$$
Notice now that 
\begin{align*}
s_{a/l}^{(2)}(\tau)&=
(2\pi\ii)^{-2}
\left(\frac {\partial^2}{\partial z^2}\right)
\log\theta(\frac al,\tau)
-
(2\pi\ii)^{-2}
\left(\frac {\partial^2}{\partial z^2}\right)_{z=0}
\log\left(\frac{\theta(z,\tau)}{z\partial_z\theta(0,\tau)}\right)\\
&=
(2\pi\ii)^{-2}\left(\frac {\partial^2}{\partial z^2}\right)
\log\theta(\frac al,\tau) - 
\frac1{12}+2\sum_{d>0}q^d\sum_{k|d}k\\
&=\frac{\ee^{2\pi\ii a/l}}
{(\ee^{2\pi\ii a/l}-1)^2}-\frac1{12}
-\sum_{d>0}q^d\sum_{k|d}k(\ee^{2\pi\ii ka/l}+\ee^{-2\pi\ii ka/l}-2),
\end{align*}
which is an Eisenstein series.  Indeed, in the notation of 
\cite[Chapter VII]{schoeneberg}, it is 
$$\frac{l^2}{4\pi^2}\left(
G_{l,2,\left(\begin{smallmatrix}0\\a\end{smallmatrix}\right)}
-G_{l,2,\left(\begin{smallmatrix}0\\0\end{smallmatrix}\right)}
\right).$$

To obtain a nice formula for $s_{a/l}^2$, recall that 
$\theta$ satisfies the heat equation
$$\theta_{zz}
=(4\pi\ii)\theta_{\tau},$$
which implies 
$$(2\pi\ii)^{-2}\left(\frac{\theta_{zz}}{\theta}\right)_z=(\pi\ii)^{-1}
\left(\frac{\theta_z}{\theta}\right)_\tau=
-4\pi\ii\sum_{d>0}dq^d\sum_{k|d}(\ee^{2\pi\ii kz}-\ee^{-2\pi\ii kz}).$$
We can integrate it with respect to $z$ while keeping in mind that
$$(2\pi\ii)^{-2}\lim_{z\to 0}
\frac{\theta_{zz}(z,\tau)}{\theta(z,\tau)}=
(2\pi\ii)^{-2}
\frac{\theta_{zzz}(0,\tau)}{\theta_{z}(0,\tau)}=
\frac14 - 6\sum_{d>0}q^d\sum_{k|d}\frac dk
$$
to obtain
$$s_{a/l}^{(2)}+s_{a/l}^2 =
(2\pi\ii)^{-2}(\frac{\theta_{zz}(a/l,\tau)}{\theta(a/l,\tau)}
-\frac{\theta_{zzz}(0,\tau)}{3\theta_{z}(0,\tau)})=
\frac 16 - 2\sum_{d>0}q^d\sum_{k|d}\frac dk(\ee^{2\pi\ii ak/l}+
\ee^{-2\pi\ii ak/l}),
$$
which is again an Eisenstein series. 
In the notation of \cite[Chapter VII]{schoeneberg}
it is equal to 
$$\frac l{2\pi^2}\sum_{a_1,a_2=0}^{l-1}
\ee^{2\pi\ii aa_1/l}
G_{l,2,\left(\begin{smallmatrix}a_1\\a_2\end{smallmatrix}\right)}.
$$
\end{proof}

\begin{proposition}\label{E2brutus}
For every even function $\chi\colon\ZZ/l\ZZ\to\CC$ the series
$$\sum_{d>0}q^d\sum_{k|d} k \,\chi(k)
\hspace{30pt}
{\rm and}
\hspace{30pt}
\sum_{d>0}q^d\sum_{k|d} \frac dk \,\chi(k)$$
lie in the linear span of $~1$, toric Eisenstein series and the
standard Eisenstein series $E_2$.
\end{proposition}

\begin{proof}
Notice that the functions $\chi_a$ defined by $k\mapsto \ee^{2\pi\ii ka}$
span the space of all functions from $\ZZ/l\ZZ$ to $\CC$, so the functions 
$\chi_a+\chi_{-a}$ span the space of all even functions from 
$\ZZ/l\ZZ$ to $\CC$. Then use explicit formulas for $s_{a/l}^{(2)}$
and $s_{a/l}^{(2)}+s_{a/l}^2$ from the proof of the above proposition.
\end{proof}

\subsection{}
Recall that $\MMM (l)$ is the space of all holomorphic modular forms for
$\Gamma_{1} (l)$ of weight two, and denote by $\EEE (l)$ the subspace of
Eisenstein series.
Proposition \ref{toricEisenstein} allows us to define a map from the 
space of Manin symbols to toric modular forms modulo Eisenstein series.
\begin{definition}
Let $M_-$ be the space of minus Manin symbols
(Definition~\ref{plusnminus}). We define a map $\mu\colon
M_-\rightarrow \MMM(l)/\EEE(l)$ by
$$
\mu((a,b)_-) = s_{a/l}s_{b/l} \mod \EEE(l)
$$
if $a,b\neq 0\mod l$, and by $0$ otherwise.
\label{toricmap}
\end{definition}

\begin{remark}
The map $\mu$ is well-defined. Indeed, it is clear that 
$(a,b)_-+(-a,b)_-$ and  $(a,b)_--(b,a)_{-}$ map to zero, even before we
mod out by $\EEE(l)$. It remains to show that
$\mu((a,b)_-+(b,c)_-+(c,a)_-)=0$ for $a+b+c=0\mod l$.  This
follows from Propositions \ref{relation} and
\ref{toricEisenstein} if $a,b,c$ are non-zero, and from Proposition
\ref{toricEisenstein} if one of them is zero.
\end{remark}

\subsection{}
We will now show that map $\mu$ commutes with the action of the Hecke 
operators $T_p$ for primes $p$ not dividing $l$. To do this,
we recall the description of toric forms in terms of cohomology
rings of smooth toric varieties \cite{BorGunn}.

Let $\deg$ be a degree function with respect to $\Sigma$ in $N$. Let
$\widehat\Sigma$ be a subdivision of $\Sigma$ such that the
corresponding toric variety $\widehat X$ is smooth.  Denote by $d_i$
the minimum nonzero lattice points of the one-dimensional cones of
$\widehat\Sigma$.  Since $\widehat\Sigma $ is a subdivision of
$\Sigma$, it is possible that some $\alpha_i=\deg(d_i)$ are now
integral.  To circumvent this difficulty, in \cite{BorGunn} we
introduced a generic degree function $\deg_1$, linear on all cones of
the original fan $\Sigma$ and such that the values
$\beta_i=\deg_1(d_i)$ are non-zero.
\begin{proposition}{\rm (\cite{BorGunn})}
In the above notation, we have
\begin{equation}
f_{N,\deg}(q)=
\lim_{\varepsilon\to 0}
\int_{{\widehat X}} \prod_i
\frac
{(D_i/2\pi\ii)\theta(D_i/2\pi\ii-\deg(d_i)-\varepsilon \deg_1(d_i),\tau
)\theta'(0,\tau )
}
{
\theta(D_i/2\pi\ii,\tau )\theta(-\deg(d_i)-\varepsilon \deg_1(d_i),\tau )
}
\label{bigugly}
\end{equation}
where $D_i$ is the cohomology class of the toric divisor corresponding
to the one-dimensional cone $\QQ_{\geq 0}d_i$ (cf. \cite{Fulton}) and
the integral means pairing with the fundamental class of $\widehat X$.
\end{proposition}

If we work modulo Eisenstein series, the above formula greatly simplifies.
\begin{proposition}
Let $d_1,\dots,d_k$ be generators of one-dimensional cones of $\widehat
\Sigma$. We denote $d_0=d_k$. Then
$$
f_{N,\deg} = \sum_{\substack
{0\leq i\leq
k-1\\\alpha_i,\alpha_{i+1}\notin\ZZ}} s_{\alpha_i}s_{\alpha_{i+1}}
\mod \EEE(l).
$$
\label{smoothtoric}
\end{proposition}

\begin{proof}
Let us first explain what happens in the case where all 
$\alpha_i$ are non-integral. We are integrating over $\widehat X$
the product of cohomology elements
$$
(1-s_{\alpha_i}D_i+r_i D_i^2),
$$
where $r_i$ is easily seen to be an Eisenstein series due to Proposition
\ref{toricEisenstein}. The intersection number $\int_{\widehat X} D_iD_j$
for $i\neq j$ equals $1$ if $d_i$ and $d_j$ come from adjacent cones
and equals $0$ otherwise. This finishes the argument.

When some of the $\alpha_i$ are zero the argument is a bit more
complicated and involves the expansion of the right hand side of 
(\ref{bigugly}) in powers of $\varepsilon$. Then up to Eisenstein series
one ends up with the integral of the product of
$(1-s_{\alpha_i}D_i)$ over $i$ for which $\alpha_i\notin\ZZ$. Details
are left to the reader.
\end{proof}

\begin{theorem}
The map $\mu\colon 
M_-\to \MMM(l)/\EEE(l)$ defined above is invariant under the 
action of the Hecke operators $T_p$ for primes $p$ coprime to $l$.
\label{muisHecke}
\end{theorem}

\begin{proof}
We will identify $N$ with the lattice of integer row
vectors with two components. Consider the fan $\Sigma$ that has
$d_1=(1,0)$, $d_2=(0,1)$, $d_3=(-1,0)$, and $d_4=(0,-1)$ as generators of
its one-dimensional faces, and a degree function on $\Sigma $
defined by $\deg(1,0)=\deg(-1,0)=m/l$ and $\deg(0,1)=\deg(0,-1)=n/l$. Proposition \ref{smoothtoric} shows that
$$f_{N,\deg}=4s_{m/l}s_{n/l} \mod \EEE(l).$$
Let us calculate the action of the Hecke operator $T_p$ on
this form modulo $\EEE(l)$. By Proposition \ref{explicitHecke}
we have
$$\stroke{f_{N,\deg}}{T_p} =
\sum_{N\subset S\subset \frac 1p N} f_{S,p\deg}.$$
Let us investigate the contribution of each $S$. To get into
the setup of Proposition \ref{smoothtoric}, we need a subdivision
$\widehat\Sigma_S$ of the fan $\Sigma$ so that the consecutive $d_i$ 
form a basis. There is a standard way of doing so. For each quadrant 
(i.e. for each two-dimensional cone of $\Sigma$), consider the set $A$ of
all non-zero points of $S$ in that quadrant. The boundary of the convex
hull of $A$ consists of two half-lines and some segments (Figure
\ref{hecke.fig}). We ignore the half lines and add to the list of $d_i$
all points in $S$ that lie on the rest of the boundary. It is easy to show
that this choice of $d_i$ guarantees that the new toric variety is
smooth. Indeed, if $d_i$ and $d_{i+1}$ did not form a basis of $S$ then
there would exist a point in $S$ lying in the convex hull of $0$, $d_i$
and $d_{i+1}$ by Pick's Theorem \cite[page 113]{Fulton}.

\begin{figure}[tbh]
\begin{center}
\includegraphics[scale = .3]{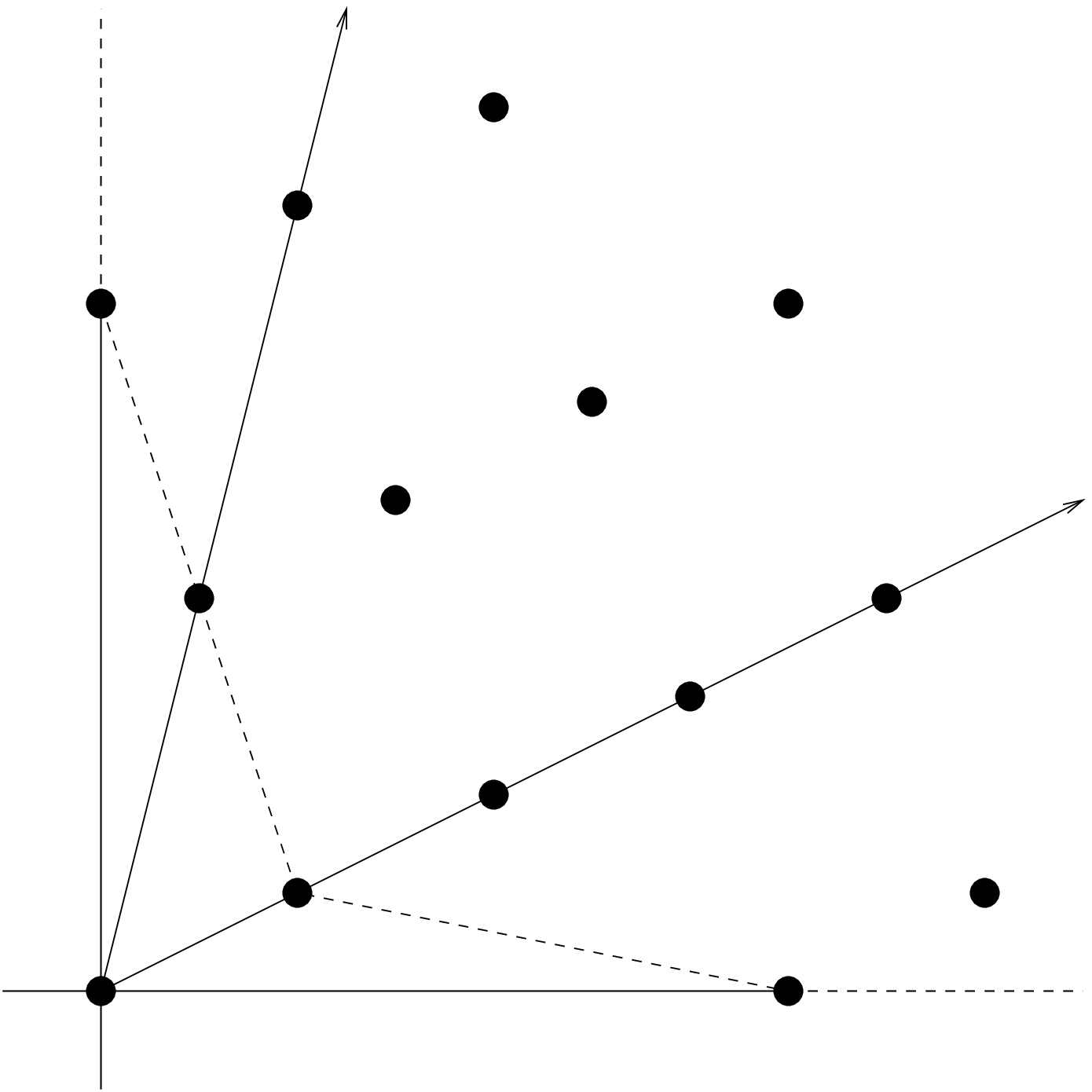}
\end{center}
\caption{\label{hecke.fig}}
\end{figure}

Because of the symmetry, it is enough to consider the 
first quadrant. Notice that if $d_i=(\frac ap,\frac cp)$ and
$d_{i+1}=(\frac bp,\frac dp)$, then $ad - bc = p$, with  $a>b\geq 0$ and
$d>c\geq 0$. Conversely, every pair of $\{(\frac ap,\frac cp),(\frac
bp,\frac dp)\}$ with $(a,b,c,d)$ as above generates some superlattice $S$ 
of coindex $p$. Moreover, this pair will form a segment of the boundary of 
the convex hull of all non-zero points of $S$ in the first quadrant.
Indeed, if any other non-zero point $(x,y)\in S$ with $x,y\geq 0$ did lie
below the line through these two points, then the area of the triangle
with vertices $(x,y)$, $(\frac ap,\frac cp)$ and $(\frac bp,\frac dp)$
would be positive but less than $\frac 1p$.

The values of $p\deg$ on the points $(\frac ap,\frac cp)$ and $(\frac
bp,\frac dp)$  are $\frac {am+cn}l$ and 
$\frac{bm+dn}l$
respectively, hence by Proposition \ref{smoothtoric} we have 
$$\stroke{f_{N,\deg}}{T_p} =
4\sum_{\substack{ad-bc=p,a>b\geq0,d>c\geq0,\\
(am+cn)/l\notin\ZZ,(bm+dn)/l\notin\ZZ}} s_{(am+cn)/l}s_{(bm+dn)/l}
\mod\EEE(l).$$
Therefore, from the definition of $\mu$, Theorem \ref{HeckeisHecke} 
and Proposition
\ref{HeckeonSymbols} we conclude that
$$\mu(T_p(m,n)_-)=\stroke{\mu((m,n)_-)}{T_p}.$$
\end{proof}

\begin{proposition}
The image of $M_-$ under $\mu$ coincides with the image of $S_-$
under $\mu$.
\label{cuspisenough}
\end{proposition}

\begin{proof}
This follows from Theorem \ref{muisHecke} and the fact that $M$ 
is isomorphic as a Hecke-module to $S_-\oplus S_+\oplus 
\EEE(l)$ \cite[Section 4.2]{Merel}.
\end{proof}

\subsection{}
Finally, we explicitly describe the composition of $\mu$ 
and the Fricke involution.
\begin{proposition}
The composition of $\mu$ and the Fricke involution $W_l$ is given by
$$
W_l\circ\mu((m,n)_-)=\tilde s_{m/l}(q)\tilde s_{n/l}(q)\mod \EEE(l)
$$
where $\tilde s_{0/l}=0$ and
$$\tilde s_{a/l}(q)=
(\frac al-\frac12 ) - \sum _{d\geq 1}q^{d}\sum _{k|d}
(\delta
_{k}^{a\mod l} - \delta
_{k}^{-a\mod l})$$
for $a=1,\ldots,l-1$.
\end{proposition}

\begin{proof}
All we need to do is calculate the Fricke involute of $s_{a/l}$.
This was accomplished in \cite{BorGunn}, up to the constant term.  To
compute the constant term, we remark that these forms have already
appeared in the literature as the \emph{Hecke-Eisenstein forms}
\cite[Chapter 15]{Lang}.
\end{proof}

\section{Forms of analytic rank zero}\label{rank.zero.section}
\subsection{}
In this subsection we define a linear map $\rho\colon \SSS(l)\to\SSS(l)$
whose image is spanned by Hecke eigenforms of analytic rank zero.
\begin{definition}
Let $f\in \SSS(l)$ be a cusp form. We define a linear map $\rho$
by
$$\rho(f)=\sum_{n=1}^{\infty}\left(\int_{0}^{\ii\infty}(\stroke{f}{
T_n})(s)ds\right) q^n.$$

\begin{proposition}\label{cuspform}
The form $\rho(f)$ is a cusp form with nebentypus equal to that of $f$.
\end{proposition}

\begin{proof}
This follows from \cite[Theorem 6]{Merel}, since $\rho (f)$ is associated to the linear map on the
Hecke algebra that maps $T$ to
$\int_{0}^{\ii\infty}(\stroke{f}{T})(s)\,ds$.
\end{proof}
\label{projection}
\end{definition}

\begin{remark}
It is easy to see that for a Hecke eigenform $f$ that is a newform
$$\rho(f)=L(f,1)f.$$
In particular, the image of the space of newforms is the span of the
new Hecke eigenforms that have analytic rank zero.
\label{newform}
\end{remark}

\begin{proposition}
The image of $\rho$ is contained in the span of all lifts of all new Hecke
eigenforms of analytic rank zero for all levels $k$, $k|l$. 
\label{imageofrho}
\end{proposition}

\begin{proof}
It is clear that if $f$ is a lift of an eigenform $g$ of
analytic rank one or more, then $\rho(f)=0$. Indeed, for every $n$
the form
$\stroke{f}{T_n}$ is a linear combination of various lifts of $g$,
which 
implies $L(\stroke{f}{T_n},1)=0$. It remains to show that if $f$ is a
lift
of $g$ then $\rho(f)$ is a linear combination of lifts of $g$.
This follows from the commutation relation
$\stroke{\rho(f)}{T_p}=\rho(\stroke{f}{T_p})$ for all prime $p$
coprime to $l$.
Indeed, from the definition of $T_p$ \cite{Lang}, 
\begin{align*}
\stroke{\rho(f)}{T_p}&=\sum_{m>0}L(\stroke{f}{  T_{mp}},1)q^m +
p\epsilon_p(\rho(f))\sum_{m>0}L(\stroke{f}{T_m},1)q^{mp}\\
&=
\sum_{m>0,(m,p)=1}L(\stroke{f}{T_{mp}},1)q^m+
\sum_{m>0, (m,p)=p}L(\stroke{f}{(T_{mp}}+p\epsilon_p(f)T_{m/p}),1)q^m\\
&=\sum_{m>0}L(\stroke{f}{T_p}T_m)q^m=\rho(\stroke f{T_p}).
\end{align*}
\end{proof}

\begin{remark}
Even though we suspect that the image of $\rho$ coincides with the 
above span, we only need the inclusion proved in the above proposition.
\end{remark}

\subsection{}
We will now prove a key result that relates the map $\rho$
to the map $W_l\circ\mu$ constructed in the previous section.
\begin{theorem}
The composition 
\begin{equation}\label{above}
\xymatrix@1{{\SSS (l)}\ar[r]^{\int }&M_{+}^{*}\ar[r]^{\pi
}&M_{-}\ar[r]^-{W_{l}\circ \mu }&\MMM (l)/\EEE (l)\ar[r]^-{\sim}&\SSS (l)},
\end{equation}
of the map induced by integration between cusps, the
map $\pi$ from
\S \ref{pi.explicit}, the map $W_l\circ\mu$, and the Hecke-equivariant
isomorphism between $\MMM(l)/\EEE(l)$ and $\SSS(l)$,
equals $\rho$.
\label{keytheorem}
\end{theorem}

\begin{proof}
Let $f$ be a cusp form, and let $\varphi\colon M_+\to\CC$ be the corresponding
element of $M_+^*$. Then by 
Theorem \ref{HeckeisHecke} and Proposition \ref{HeckeonSymbols},
\begin{multline*}
\rho(f)(q)=
\sum_{n>0}q^n\langle \stroke{f}{T_n},(1,0)\rangle=
\sum_{n>0}q^n\langle f,T_n(1,0)\rangle\\
=-\sum_{n>0}q^n\varphi(T_n(0,1))=
-\sum_{n>0}q^n\sum_{\substack{ad-bc=n,\\a>b\geq 0,d>c\geq 0}}
\varphi((c,d)_+).
\end{multline*}
Notice that we are using our convention that $\varphi(c,d)=0$ for 
$g.c.d.(c,d,l)>1$. On the other hand, by Theorem \ref{2.9}
the composition of all the maps
in \eqref{above} except for the last one yields an element of
$\MMM(l)/\EEE(l)$
given by 
$$\rho_1(f)(q)=\frac 1{12}
\sum_{a,b=0}^{l-1}(\varphi((a,a-b)_+)-\varphi((a,a+b)_+))
\tilde s_{a/l}(q) \tilde s_{b/l}(q),$$ where we again apply our
convention.  Also, we can formally use the same expression for $\tilde
s_{0/l}$ as for the rest of $\tilde s_{a/l}$ because the coefficient
at $\tilde s_{0/l}\tilde s_{b/l}$ is zero.

In what follows, it will be convenient for us to ignore the constant terms
in all our expressions. Indeed, all our functions are modular forms, so
any constant term can always be restored. We will denote all these 
constant terms by $C$. Using 
$$\tilde s_{a/l}\tilde s_{b/l} = C + \sum_{n>0}q^n
\Bigl(
\sum_{k|n,k>0}(\frac 12-\frac al)(\delta_k^{b\mod
l}-\delta_k^{-b\mod l}) +
\Bigr.
$$
$$
\Bigl.
\sum_{k|n,k>0}(\frac 12-\frac bl)(\delta_k^{a\mod
l}-\delta_k^{-a\mod l})+
\sum_{\substack{m_1k_1+m_2k_2=n,\\m_1,k_1,m_2,k_2>0}}
(\delta_{k_1}^{a\mod l}-\delta_{k_1}^{-a\mod l})
(\delta_{k_2}^{b\mod l}-\delta_{k_2}^{-b\mod l})
\Bigr)
$$
and symmetry properties of $\varphi((a,a-b)_+)-\varphi(a,a+b)_+)$, we get
$$
\rho_1(f)(q) = C + \frac13\sum_{n>0}q^n
\Bigl(
\sum_{k|n,k>0}\sum_{a=0}^{l-1}(\frac 12 -\frac al)
(\varphi((a,k-a)_+)-\varphi((a,k+a)_+))
\Bigr.
$$
$$
\Bigl.
+\sum_{\substack{m_1k_1+m_2k_2=n,\\m_1,k_1,m_2,k_2>0}}
(\varphi((k_1,k_1-k_2)_+)-\varphi((k_1,k_1+k_2)_+))
\Bigr).
$$
Let us now simplify the second part of this expression to make it look
more like $\rho(f)$. 
We split it into four sums
\begin{equation}
\begin{split}
\sum_{\substack{m_1k_1+m_2k_2=n,\\m_1,k_1,m_2,k_2>0\\k_1\geq k_2}}
\varphi((k_1,k_1-k_2)_+)
&+
\sum_{\substack{m_1k_1+m_2k_2=n,\\m_1,k_1,m_2,k_2>0\\k_1< k_2}}
\varphi((k_1,k_1-k_2)_+)
\\
-
\sum_{\substack{m_1k_1+m_2k_2=n,\\m_1,k_1,m_2,k_2>0\\m_1> m_2}}
\varphi((k_1,k_1+k_2)_+)
&-
\sum_{\substack{m_1k_1+m_2k_2=n,\\m_1,k_1,m_2,k_2>0\\m_1\leq m_2}}
\varphi((k_1,k_1+k_2)_+)
\end{split}
\label{foursums}
\end{equation}
and deal with each sum separately. We will give a detailed calculation
for one of the sums and will indicate how to manipulate the other three.
\begin{lemma}
$$
\sum_{\substack{m_1k_1+m_2k_2=n,\\m_1,k_1,m_2,k_2>0\\m_1> m_2}}
\varphi((k_1,k_1+k_2)_+)
=
\sum_{\substack{m_1k_1+m_2k_2=n,\\m_1,k_1,m_2,k_2>0\\m_1>m_2}}
\varphi((k_1,k_2)_+)
+
\sum_{\substack{ad-bc=n,\\a>b>0,d>c>0}}
\varphi((c,d)_+)
$$
\label{sumthree}
\end{lemma}

\begin{proof}[Proof of the lemma] We first of all rewrite
$\varphi((k_1,k_1+k_2)_+)$ as $\varphi((k_1,k_2)_+) +
\varphi((k_2,k_1+k_2)_+) $. Then for the second term we make the
change of variables $(a,b,c,d)=(m_1,m_1-m_2,k_2,k_1+k_2)$.
\end{proof}

We perform similar but easier manipulations for the remaining three sums
in (\ref{foursums}). For the first sum we make the change of variables
$(a,b,c,d)=(m_1+m_2,m_2,k_1-k_2,k_1)$. For the second sum we
make the change of variables $(m_1',m_2',k_1',k_2')=
(m_1+m_2,m_2,k_1,k_2-k_1)$, so that it cancels the first sum in
Lemma \ref{sumthree}. For the fourth sum in (\ref{foursums})
we make the change of variables $(a,b,c,d)=(m_2,m_2-m_1,k_1,k_1+k_2)$.
After some straightforward calculations we get
$$
\sum_{\substack{m_1k_1+m_2k_2=n,\\m_1,k_1,m_2,k_2>0}}
(\varphi((k_1,k_1-k_2)_+)-\varphi((k_1,k_1+k_2)_+))
=
-3\sum_{\substack{ad-bc=n,\\a>b\geq0,d>c\geq0}}
\varphi((c,d)_+)
$$
$$
-\sum_{0<k|n}(\frac {2n}k+1)\varphi((k,0)_+)
-2\sum_{\substack{k_1|n,k_1>k_2>0}}
\varphi((k_1,k_1-k_2)_+)
.$$
As a result, 
$$\rho_1(f)(q)-\rho(f)(q)
=C+
\frac 13\sum_{n>0}q^n
\sum_{k|n,k>0}
\Bigl(
\sum_{a=0}^{l-1}(\frac 12 -\frac al)
(\varphi((k,a-k)_+)-\varphi((k,a+k)_+))
\Bigr.
$$
$$
\Bigl.
-(\frac {2n}k-1)\varphi((k,0)_+)
-2\sum_{0\leq m \leq k}
\varphi((k,m)_+)
\Bigr).
$$
It is possible to further simplify this equation to obtain
$$\rho_1(f)(q)-\rho(f)(q)
=C - \frac23\sum_{n>0}q^n\sum_{k|n,k>0}\frac nk \varphi((k,0)_+) -
\frac 2{3l}\sum_{n>0}q^n\sum_{k|n,k>0} k
\Bigl(\sum_{b=0}^{l-1}\varphi((k,b)_+)\Bigr),
$$
which is an Eisenstein series. Indeed, it can be easily written (up to
a constant) as a linear combination of $s_{a/l}^{(2)}$, $s_{a/l}^2$
and the (non-modular) $SL_2(\ZZ)$-Eisenstein series $E_2$ by
Proposition \ref{E2brutus}. It remains
to observe that the coefficient of $E_2$ must be zero, because of
the transformation properties of $E_2$ under $\Gamma_1(l)$.
\end{proof}

\subsection{}
We are now ready to prove our main result.
\begin{theorem}
For each integer $l$ the space spanned by pairwise products of $s_{a/l}$ for
$a=1,\dots ,l-1$ is the direct sum of the span of all Hecke eigenforms of
analytic rank zero and some subspace of the space of Eisenstein
series.  \label{main}
\end{theorem}

\begin{proof}
We will prove this theorem by induction on $l$. For small levels there
is nothing to prove, because there are simply no cusp forms.

Fix $l$ and assume that the statement of the theorem is true for all
smaller levels. In particular, this implies that lifts of all forms of
smaller levels are contained in the span of toric forms, because the
space of toric forms is stable under liftings, see Theorem
\ref{summary}.  Every new Hecke eigenform $f$ of analytic rank zero is
contained in the image of $\rho$, see Remark \ref{newform}. By
Theorem \ref{keytheorem}, $f$ is contained in the image of $W_l\circ
\mu$, and so is a toric form up to Eisenstein series. Because the
space of toric forms is Hecke stable, this implies that $f$ is
toric.  This proves that the space of toric forms contains the span of
all Hecke eigenforms of analytic rank zero.

To prove the opposite inclusion, notice that by the induction assumption
it is enough to consider $s_{a/l}s_{b/l}$ with $g.c.d.(a,b,l)=1$.
By Proposition \ref{cuspisenough}, there is an element
$x\in S_-$ such that $\mu(x)=s_{a/l}s_{b/l}$. We use here that 
$g.c.d.(a,b,l)=1$, because otherwise the symbol $(a,b)_-$ is not
defined. By the definition of $\pi$ in \S \ref{pi.explicit}
there is an element
$\varphi\in S_+^*$ such that $\pi(\varphi)=x$. Moreover, we can
find a cusp form $f$ which induces the linear map $\varphi$ on 
$S_+$. Then Theorem \ref{keytheorem} shows that $s_{a/l}s_{b/l}$
is proportional to $W_l\circ\rho(f)$ up to an Eisenstein series.
By Proposition \ref{imageofrho} $s_{a/l}s_{b/l}$ lies in the span of 
Hecke eigenforms of rank zero and Eisenstein series, which finishes
the proof. 
\end{proof}

\begin{remark}
It is easy to see that as a corollary we get similar results for any given
nebentypus. In particular, modulo Eisenstein series, the span of the
forms 
$$\sum_{k\in(\ZZ/l\ZZ)^*} s_{ka}s_{kb}, \quad a,b=1,\dots ,l-1$$ 
is the span of all Hecke eigenforms of analytic rank zero for
$\Gamma_0(l)$.
\end{remark}

\begin{remark}
In general it is not clear which Eisenstein series are toric forms.
In particular at level $l=25$ the space of toric Eisenstein series
has codimension one.
\end{remark}

%
%

\bibliographystyle{amsplain}
\bibliography{vanish}

\end{document}